\documentclass{article}
\usepackage[utf8]{inputenc}
\usepackage{amsmath}
\usepackage{amsfonts}
\usepackage{mathrsfs}

\usepackage[margin=1in]{geometry}
\usepackage{graphicx}
\usepackage[thinlines]{easytable}
\usepackage{url}
\usepackage{authblk}

\title{Blind Number Sequencing}
\author{Parker Kuklinski, Nick Vogel}
\date{October 2023}

\begin{document}

\maketitle

\begin{abstract}
A popular challenge on the social media app TikTok is to place the output of a random number generator in ascending order without arriving at a contradiction. Most players rely on intuition to construct their sequences and debate has surfaced in the comments about the optimal strategy to win the game. To this end, we determine an optimal strategy that offers a 21\% improvement in winning odds over the intuitive ``equal-spacing" strategy in the popular $\#$20numberchallenge variant. 
\end{abstract}

\section{Introduction}

TikTok is a video hosting social media platform owned by ByteDance \cite{herrman19}, and according to similarweb was the $14^\text{th}$ most popular website on the internet in October 2023 \cite{similarweb}. Through certain design choices, The platform has curated a unique culture that encourages users to follow trends, often denoted with the hashtag symbol, \#. Certain trends are facilitated through the use of `filters', or video overlays that warp faces, modify voices, or even act as interactive games. One simple filter that became popular is the random number filter, which at the tap of the screen, displays a random number between 0-999 over the user's head. Many use the filter as an oracle to answer humorous questions (e.x. How many years do I have left to live? How many push-ups should I do right now? How many kids will I have?)

One user eventually invented the 5-number challenge, played as follows: the user begins with an empty list of five items, generates a random number, and places that number somewhere in the list. The user then continues to generate random numbers one at a time and placing them in empty slots in the list \emph{in ascending order} until either the number cannot fit in any available slots or all slots have been filled. This would sometimes be referred to as `blind number sequencing' since the player is attempting to sequence numbers one at a time without knowing what comes next. The game was conceptually simple, entertaining to watch, easy to play along, and thus primed to become a popular trend in the community.

As users started beating the 5-number challenge with ease (win probability with optimal strategy is $\approx 22\%$), participants considered larger number games, eventually coalescing around the 20-number challenge. To date, only a handful of players have ever succeeded in completing the challenge. Due to the difficulty of the challenge and mass of viewers, debate over optimal strategies raged in the comments of these $\#$20numberchallenge videos; some argued that placing numbers based on evenly dividing the interval was optimal, others advocated for riskier maneuvers to remain in the game longer.

Fortunately the $n$-number challenge is a probabilistic system that can be analytically modeled; with some effort optimal strategies and exact win probabilities can be computed. In particular, the probability of winning the 20-number challenge with an optimal strategy is just under $1/8000$. It turns out that this optimal strategy is not to evenly divide the interval as some suggested, but to be slightly more avoidant of placing numbers at the ends of available bins.

The remainder of the paper is structured as follows: section 2 formally defines an $n$-number game and what is meant by a strategy. Section 3 presents the equal-spacing strategy while section 4 documents the optimal risk-tolerant strategy. Section 5 notes a practical consideration of implementing these strategies that slightly tempers the advantage of the risk-tolerant strategy.

\section{Definitions and Methods}

Let $n\in\mathbb{N}$. In the $n$-\emph{number game}, a random sequence ${\bf x}\in \{ 000,...,999\} ^n$ (not necessarily in ascending order) is revealed to the player one element at a time, who likewise defines a function $\sigma :{\bf x} \rightarrow\{ 1,...,n\}$ one element at a time. The player \emph{wins} if $\sigma$ is a monotonically increasing bijection. To ease the analysis in this paper, we consider a variant where the elements ${\bf x}\in [0,1]^n$ are uniformly distributed random variables. We say a \emph{strategy} $\sigma=\{\sigma _1,...,\sigma _n\}$ is a collection of functions $\sigma _k:[0,1]\rightarrow\{ 1,...,k\}$ that each map a random number to an entry on a length $k$ list. These strategies can be represented as nondecreasing step functions such that for each strategy function $\sigma _k$ we can associate a collection of numbers $\{\alpha _{k,0},...,\alpha _{k,k}\}$ such that if $x\in (\alpha _{k,j-1},\alpha _{k,j})$, then $\sigma _k(x)=j$ (we necessarily have $\alpha _{k,0}=0$ and $\alpha _{k,k}=1$).

The objective of the $n$-number game is to then define a strategy $\sigma$ (or equivalently a collection of $\{\alpha _{k,j}\} _{1\le j\le k\le n}$) which maximizes the probability $p_n$ of correctly blind-sequencing an $n$-length random sequence ${\bf x}$. Suppose $x\in {\bf x}$ is the first random number presented in the list. If we place this number in the $\sigma _n(x)$ slot in our empty sequence, then the probability we can correctly blind sequence the rest of the numbers in ${\bf x}$ is equal to the probability that $\sigma _n(x)$ is the correct placement of $x\in {\bf x}$ (otherwise that there are $\sigma _n(x)-1$ elements of ${\bf x}$ less than $x$ and $n-\sigma _n(x)$ elements of ${\bf x}$ greater than $x$) times the probability that we will be able to sequence those uniformly distributed random elements as they are presented to us one at a time (i.e. times $p_{\sigma _n(x)-1}p_{n-\sigma _n(x)}$). Integrating $x$ across the domain $[0,1]$ gives a recursive convolution-type definition of $p_n$ given a strategy $\tilde{\sigma}$:
$$p_n=\int _0^1 {n-1 \choose \sigma _n(x)-1}x^{\sigma _n(x)-1}(1-x)^{n-\sigma _n(x)}p_{\sigma _n(x)-1}p_{n-\sigma _n(x)} dx$$
Here we adopt convention that $p_0=1$. By recognizing the piecewise structure of $\sigma _n(x)$ we can clean up the notation as follows:
\begin{equation}
    p_n=\sum _{k=1}^n{n-1 \choose k-1}p_{k-1}p_{n-k}\int _{\alpha _{n,k-1}}^{\alpha _{n,k}} x^{k-1}(1-x)^{n-k} dx
\end{equation}
It may be tempting to attempt to collect these probabilities in a generating function $p(z)=\sum _{n=0}^\infty p_nz^n$ and attempt derive a formula for $p(z)$ using Gosper-Zeilberger algorithm \cite{zeilberger91}, however a closed form is unreachable.

\section{Equal-Spacing Strategy}

Without explicitly calculating optimal values of $\alpha _{k,j}$ (i.e. maximizing $p_n$), one may naively assume that an equal-spacing strategy would be optimal, or that $\alpha _{k,j}=j/k$. One can show that this strategy is greedy in the sense that it maximizes the probability that a sequence in progress is correct \emph{so far}. Let $v_n\in \{ [0,1]\cup\emptyset\} ^n$ denote an $n$-length sequence in progress where if $v_n$ has $k$ nonempty entries with corresponding inputs in ${\bf x}$, this represents the current state of an $n$-number game at the $k^\text{th}$ turn. We can then define a function $\tilde{p}:\{ [0,1]\cup\emptyset\} ^n\rightarrow [0,1]$ such that $\tilde{p}(v_n)$ is the probability that $v_n$ correctly orders its nonempty elements in ${\bf x}$. If the nonempty elements of $v_n$ partition $\{ 1,...,n\}$ into $j$ "bins" of consecutive empty elements, let $n_i$ be the size of the $i^\text{th}$ bin and let $p_i$ be the probability that a random element $x\in [0,1]$ will fit into the $i^\text{th}$ bin. Then $\tilde{p}(v_n)$ is simply the probability that the remaining $n-k$ elements of ${\bf x}$ can be placed in the $j$ empty bins, otherwise
$$\tilde{p}(v_n)={n_1+\hdots +n_j \choose n_1,...,n_j}p_1^{n_1}\cdot\hdots\cdot p_j^{n_j}$$

We can prove that an equal-spacing strategy $\sigma ^{ES}$ maximizes $\tilde{p}$ for every new element $x$ appended to $v_n$. Due to the product structure of $\tilde{p}(v_n)$ and since each new $x$ can only fit in one of the $j$ empty bins, it suffices to prove optimality for the fully empty vector $v_n=\emptyset$. If we are presented with a new random element $x\in [0,1]$, let $v_n(x)$ be the updated vector with $x$ in the $\sigma _n(x)=k$ slot. Then we have
$$\tilde{p}(v_n(x))={n-1\choose k-1}x^{k-1}(1-x)^{n-k}$$
For shorthand, let $\tilde{p}_k$ be this probability corresponding to selecting $\sigma _n(x)=k$. Then we have
$$\frac{\tilde{p}_{k+1}}{\tilde{p}_k}=\frac{(n-k)x}{k(1-x)}$$
It follows that $\tilde{p}_{k+1}>\tilde{p}_k$ when $x>k/n$ and $\tilde{p}_k>\tilde{p}_{k+1}$ when $x<k/n$. Therefore $\tilde{p}_k$ is maximized when we choose $x\in (\frac{k-1}{n},\frac{k}{n})$.

\section{Risk-Tolerant Strategy}

While the equal-spacing strategy aligns with intuition and greedily maximizes the probability that the current vector is correct, it does not take into account the probability that one will be able sequence the remaining elements as they are presented. In this section, we will demonstrate that in certain cases, selecting a less likely placement of $x\in {\bf x}$ is countered by an improvement in the probability of being able to order the remaining elements of ${\bf x}$. Adopting these modified strategies results in small gains in $p_n$ for $n$ small, but significant gains as $n$ increases.

Let us first demonstrate a superior strategy to equal-spacing in the case of $n=3$. It is easy to show that a symmetric selection criteria for $n=2$ (i.e. $\alpha _{2,1}=1/2$) is still optimal with $p_2=3/4$. For $n=3$ let $\alpha _{3,1}=1-\alpha _{3,2}=\alpha$. It can be shown that for an equal-spacing strategy $\alpha =1/3$ we have $p_3=83/162\approx 0.5123$. However revisiting (1) we can rewrite $p_3$ in terms of variable $\alpha$:
$$p_3=\frac{11}{6}a^3-\frac{7}{2}a^2+\frac{3}{2}a+\frac{1}{3}$$
From this formula we see that $p_3$ attains a maximum value of $p_3=377/726\approx 0.5193$ when $\alpha =3/11$. Intuitively speaking, rather than maximizing the probability that the random number $x$ belongs in slot $k$, we sacrifice some of that probability so that if $x$ is indeed in slot $k$, sequencing the rest of the numbers will be easier. In the case of the 3-number game, we are more inclined to place the first random number $x$ in the $2^{nd}$ window because, if correct, the remaining numbers are automatically sequenced depending on whether they lie below or above $x$. This improvement in odds of winning over the equal-spacing strategy is relatively small for the 3-number game, however we will see that the improvement becomes more significant as $n$ increases.

Fortunately, there is a readily apparent visual intuition to selecting the optimal strategy for larger $n$. Revisiting EQN, let us define the function
\begin{equation}
    f_{n,k}(x)={n-1 \choose k-1}p_{k-1}p_{n-k}x^{k-1}(1-x)^{n-k}
\end{equation}
or the probability of correctly sequencing $n$ uniform random variables on $[0,1]$ given that you have placed the first random variable $x$ in the $k^\text{th}$ slot. From the piecewise structure of (1), it is evident that in order to maximize $p_n$, we need to select $\{ a_{n,k}\}$ such that for $x\in (\alpha _{n,k-1},\alpha _{n,k})$, $f_{n,k}(x)>f_{n,j}(x)$ for all $j\ne k$. This is depicted in Figure 1 as we plot all of the functions $f_{n,k}(x)$ for fixed $n$. It thus follows that the optimal $\alpha _{n,k}$ are satisfied by the equations $f_{n,k}(\alpha _{n,k})=f_{n,k+1}(\alpha _{n,k})$. This affords us the following closed form:
\begin{equation}
    \alpha _{n,k}=1/\left(1+\frac{p_kp_{n-k-1}}{p_{k-1}p_{n-k}}\left(\frac{n}{k}-1\right)\right)
\end{equation}

\begin{figure}
\centering
\includegraphics[scale=0.37]{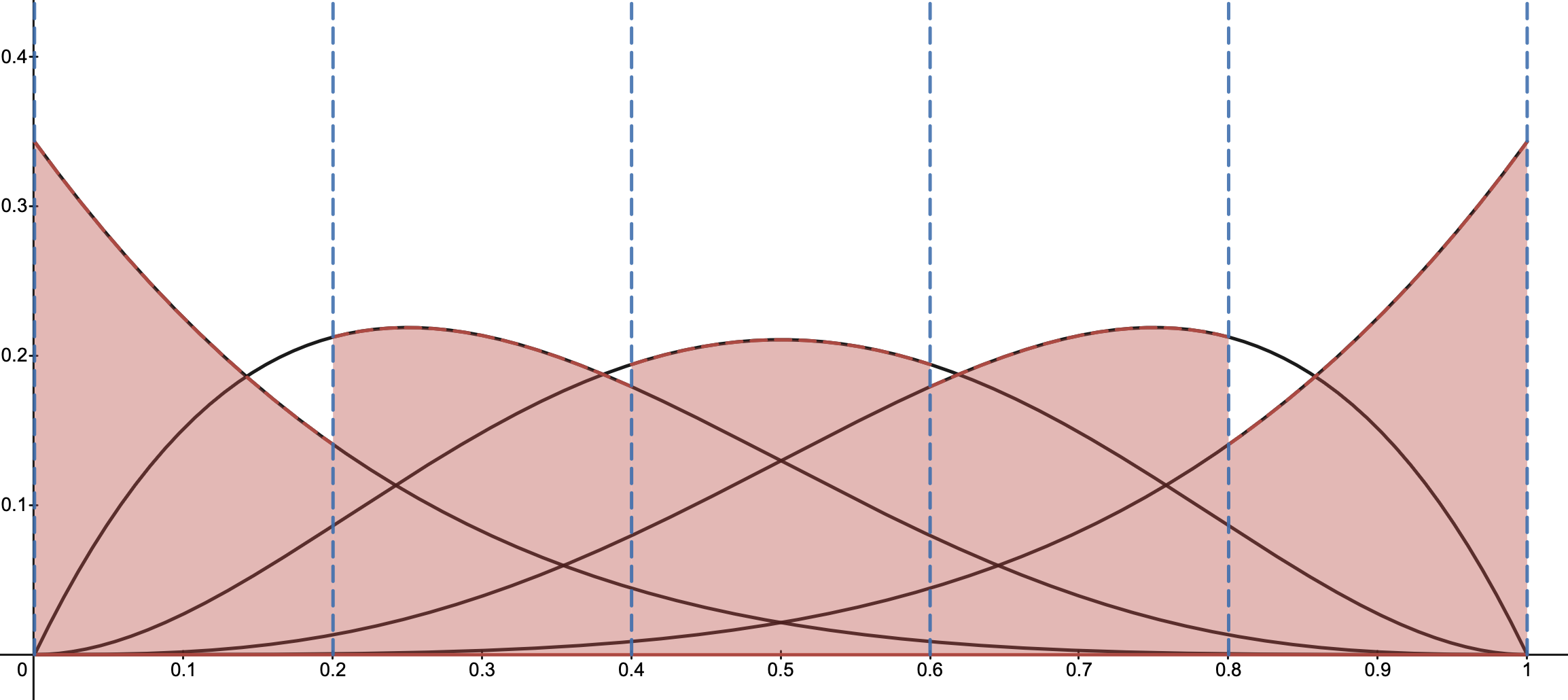}

\vspace{1cm}

\includegraphics[scale=0.37]{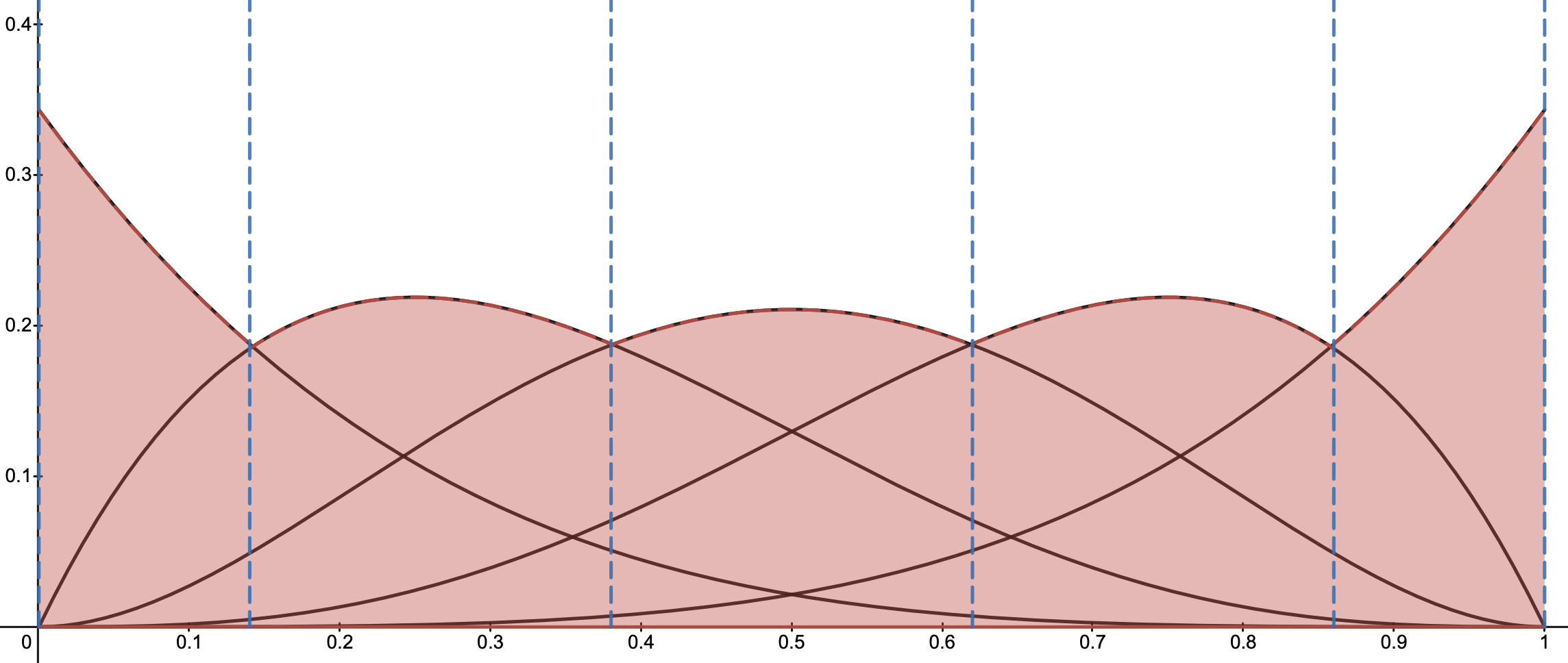}
\caption{Illustration of advantage of risk-tolerant strategy over equal-spacing strategy for a 5-number game. Black curves correspond to $f_{5,k}(x)$ and blue vertical dashed lines indicate strategy boundaries. Volume of red area equals probability of winning the 5-number game given the blue decision boundaries and optimal play after. (\emph{Top}) Equal-spacing strategy (\emph{Bottom}) Risk-tolerant strategy. Notice entire volume filled in the bottom figure.}
\end{figure}

\begin{table}
\begin{center}
\begin{tabular}{ |p{1cm}|p{1cm}||p{2cm}|p{2cm}|p{2cm}|  }
\hline
 Order & Sequence so far & Prob correct  so far& Win prob (RT) & Win prob (ES) \\
 \hline
 1 & ~ & ~ & ~ & ~ \\
 \hline
2 & ~ & ~ & ~ & ~ \\
 \hline
 3 & 130 & ~ & ~ & ~ \\
 \hline
 4 & ~ & $9.58\times 10^{-3}$ & $1.12\times 10^{-4}$ & $1.02\times 10^{-4}$ \\
 \hline
 5 & ~ & $6.81\times 10^{-3}$ & $1.28\times 10^{-4}$ & $1.17\times 10^{-4}$ \\
 \hline
 6 & ~ & $2.07\times 10^{-3}$ & $4.61\times 10^{-5}$ & $4.26\times 10^{-5}$ \\
 \hline
 7 & ~ & $3.51\times 10^{-4}$ & $8.37\times 10^{-6}$ & $7.74\times 10^{-6}$ \\
 \hline
 8 & ~ & $3.56\times 10^{-5}$ & $8.50\times 10^{-7}$ & $7.85\times 10^{-7}$ \\
 \hline
 9 & ~ & $2.17\times 10^{-6}$ & $4.81\times 10^{-8}$ & $4.46\times 10^{-8}$ \\
 \hline
 10 & ~ & $7.33\times 10^{-8}$ & $1.37\times 10^{-9}$ & $1.25\times 10^{-9}$ \\
 \hline
 11 & ~ & $1.06\times 10^{-9}$ & $1.25\times 10^{-11}$ & $1.13\times 10^{-11}$ \\
 \hline
 12 & 573 & ~ & ~ & ~ \\
 \hline
 13 & ~ & ~ & ~ & ~ \\
 \hline
 14 & ~ & ~ & ~ & ~ \\
 \hline
 15 & ~ & ~ & ~ & ~ \\
 \hline
 16 & 761 & ~ & ~ & ~ \\
 \hline
 17 & ~ & ~ & ~ & ~ \\
 \hline
 18 & ~ & ~ & ~ & ~ \\
 \hline
 19 & ~ & ~ & ~ & ~ \\
 \hline
 20 & ~ & ~ & ~ & ~ \\
 \hline
\end{tabular}
\caption{Illustration of a 20-number challenge where probabilities displayed correspond to placing the next number 170 in slot $k$. Notice that the equal-spacing strategy (i.e. deciding based on maximizing probability correct so far) would dictate that 170 should be placed in slot 4, whereas maximizing probability of winning using the risk-tolerant strategy suggests that 170 should be placed in slot 5. By placing 170 in slot 5 and continuing to play risk-tolerant, the odds of winning are about 25\% greater than placing 170 in slot 4 and continuing to play with equal-spacing.}
\end{center}
\end{table}

We present a variety of illustrations to elucidate the advantage of using the risk-tolerant strategy versus the equal-spacing strategy. Let $p_n^{ES}$ and $p_n^{RT}$ correspond to win probabilities of the $n$-number game using the equal-spacing and risk-tolerant strategies respectively. Figure 2 demonstrates that both $p_n^{ES}$ and $p_n^{RT}$ approximately exponentially decrease as $n$ increases, and further the factor of advantage of risk-tolerant over equal-spacing (i.e. $p_n^{RT}/p_n^{ES}$) is approximately linear in $n$. For example, in a 20-number game the equal spacing strategy gives a win probability of $p_{20}^{ES}\approx 1/9,651$ while the risk-tolerant strategy achieves $p_{20}^{RT}\approx 1/7,980$, the improvement factor being $p_{20}^{RT}/p_{20}^{ES}\approx 1.2095$. The advantage is more pronounced for larger number games; in a 40-number game the win probabilities satisfy $p_{40}^{ES}\approx 1/433,000,000$, $p_{40}^{RT}\approx 1/296,000,000$, and an improvement factor of $p_{40}^{RT}/p_{40}^{ES}\approx 1.4617$

\begin{figure}
\centering
\includegraphics[scale=0.5]{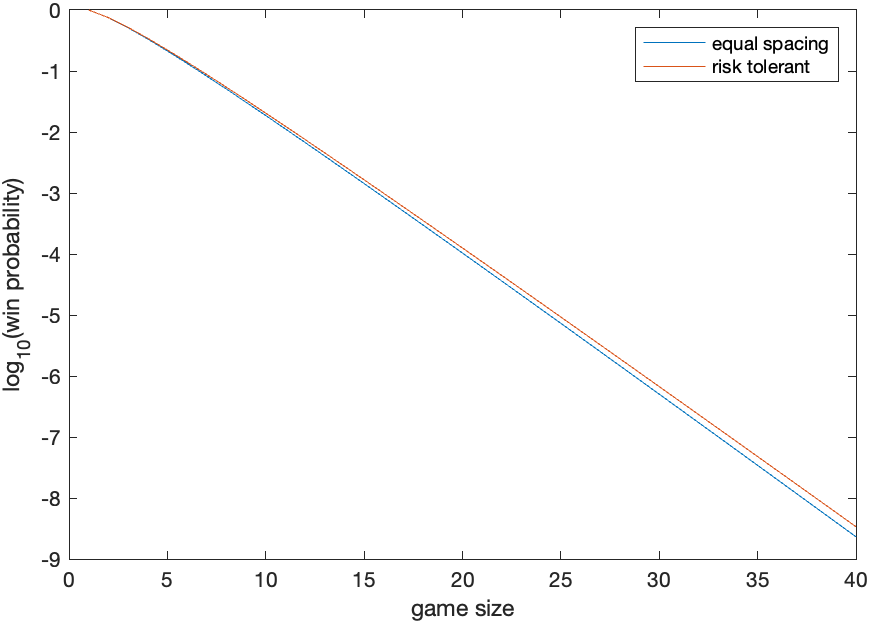}\hspace{1cm}\includegraphics[scale=0.5]{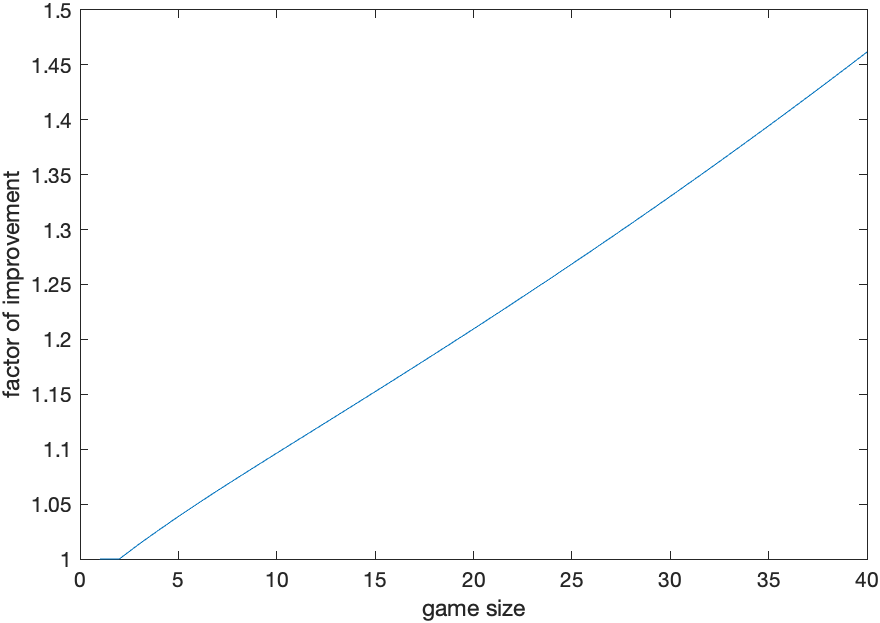}
\caption{Comparison of equal-spacing strategy versus risk-tolerant strategy (\emph{Left}) log-scale probability of winning (\emph{Right}) factor of improvement in win probability of risk-tolerant over equal-spacing.}
\end{figure}

Figures 3 and 4 illustrate the decision boundaries corresponding to each strategy, which we denote as $\alpha _{n,k}^{ES}$ and $\alpha _{n,k}^{RT}$. In particular, we find that the most extreme differences in bin sizes between equal-spacing and risk-tolerant strategies occur in the first and last bins, namely the size of the end bins in the risk tolerant strategy is about 60\% of the size of an equal-spacing bin (i.e. $\lim _{n\rightarrow\infty}n\alpha _{n,1}\approx 0.6$). Meanwhile, the second and third risk-tolerant bins (equivalently second and third to last bins) have volume slightly less than $1/n$ each, and the remainder of the bins approximately evenly divide the remainder.

\begin{figure}
\centering
\includegraphics[scale=0.42]{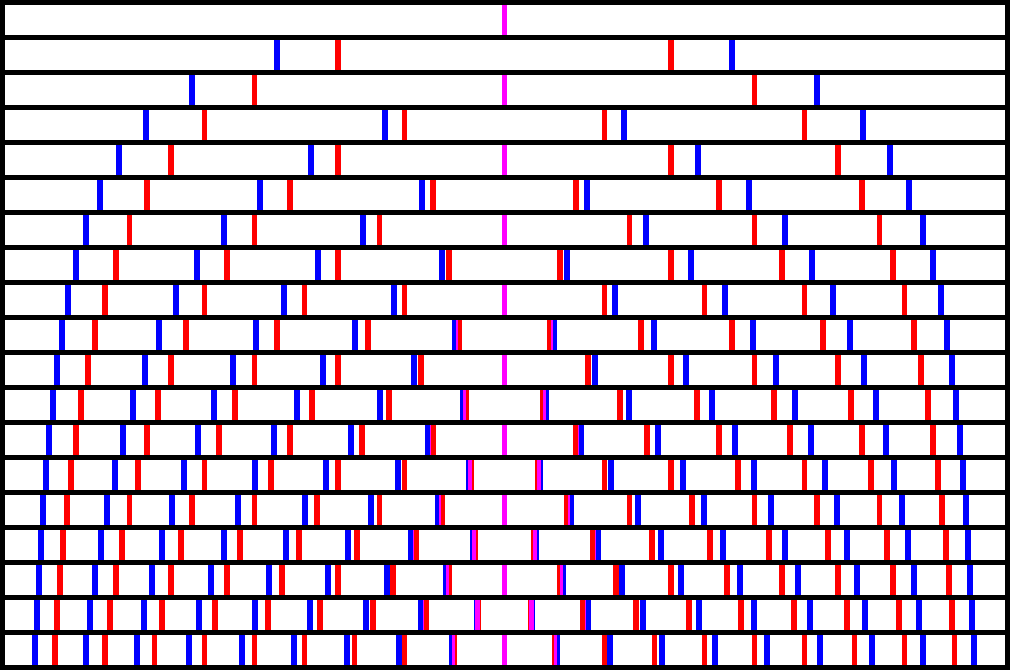}
\caption{Comparison of bin boundaries for different strategies. Each row indicates an $n$-number game where $n\in\{ 2,...,20\}$. Blue ticks indicate boundary for risk-tolerant strategy, red ticks indicate boundary for equal spacing strategy.}
\end{figure}

\begin{figure}
\centering
\includegraphics[scale=0.7]{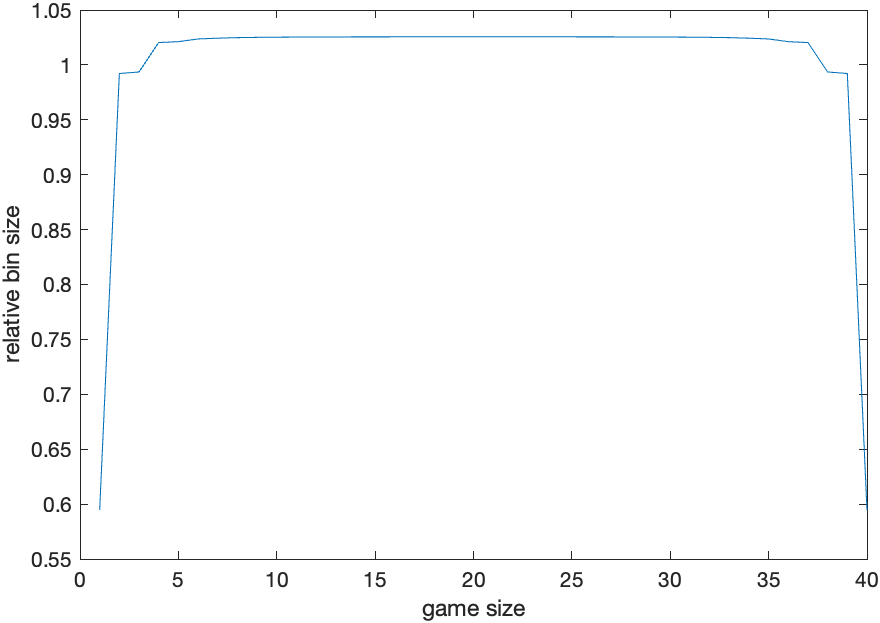}
\caption{Relative bin size (as a factor of the equal-spacing bin size $1/n$) of the risk-tolerant strategy for $n=40$.}
\end{figure}

While the equal-spacing strategy may seem more intuitive, the risk-tolerant strategy suggests that one should exercise caution when placing numbers at the ends of the distribution. Interestingly, the second author has played tens of thousands of 20-number games, and prior to any formal mathematical analysis has developed their own instinctive strategy strongly mimicking the risk-tolerant strategy.

\section{Considerations for multiple games}

We have just demonstrated that the risk-tolerant strategy outperforms the equal-spacing strategy in likelihood of winning, but we are missing a key practical component of the game. Since winning the 20-number challenge is exceedingly rare, one would likely play until they win once, and since this is a time-intensive endeavor, more important than maximizing the win probability of a single game is \emph{minimizing} the expected time playing until the game is won. Let $f_n^{\sigma}:\{ 1,...,n\}\rightarrow [0,1]$ be a conditional probability distribution where $f_n^{\sigma}(k)$ is the probability that the player using strategy $\sigma$ is eliminated at the $k^\text{th}$ turn \emph{given} that the player does not win. Then from geometric series arguments it follows that $\mathbb{E}_n(\sigma)$, the expected number of random elements drawn using strategy $\sigma$ before winning a single $n$-number challenge is given by:
$$\mathbb{E}_n(\sigma)=n+\left(\frac{1}{p_n^\sigma }+1\right)\mathbb{E}(f_n^\sigma )$$

While increasing the probability of success through a more optimal strategy will intuitively decrease the expected time playing the $n$-number challenge before reaching a first victory, we must be careful that this is not offset by an increase in $\mathbb{E}(f_n^\sigma )$, or otherwise that we are not spending too much time playing rounds that do not eventually win. We can see in Figure 5 that while $p_n^{RT}>p_n^{ES}$, the graph of $f_{n}^{RT}$ is shifted to the right of $f_{n}^{ES}$. In particular, we estimate $\mathbb{E}\left( f_{20}^{ES}\right) \approx 10.8414$ and $\mathbb{E}\left( f_{20}^{RT}\right) \approx 11.5178$, so that even though an individual 20-number game using a risk-tolerant strategy is about 21\% more likely to win than an equal-spacing strategy, a risk-tolerant player will on average play only 13\% fewer number draws before a first win than an equal-spacing player will. A similar calculation shows that for a 40-number challenge, the risk-tolerant player uses on average 26\% fewer number draws before winning than the equal-spacing player.

\begin{figure}
\centering
\includegraphics[scale=0.5]{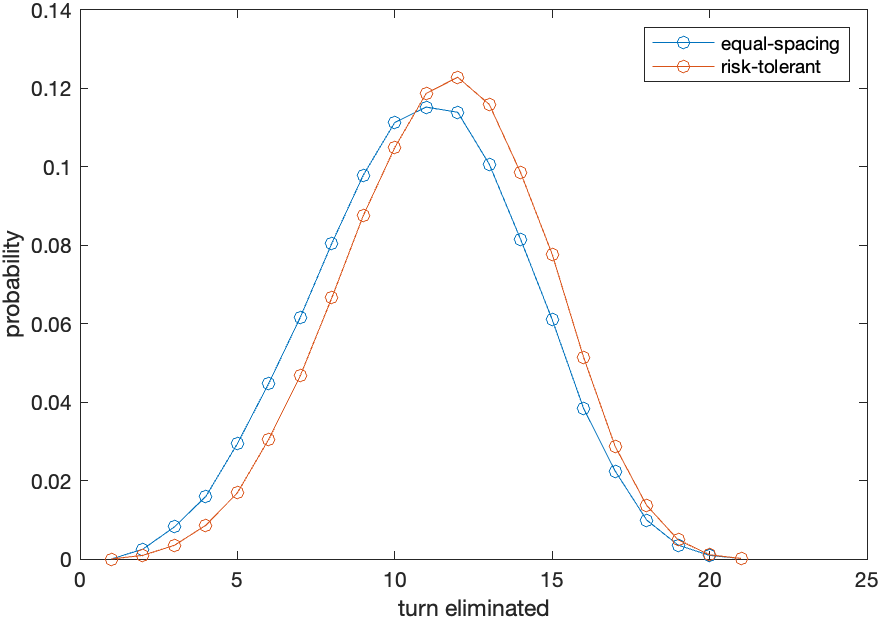}\hspace{1cm}\includegraphics[scale=0.5]{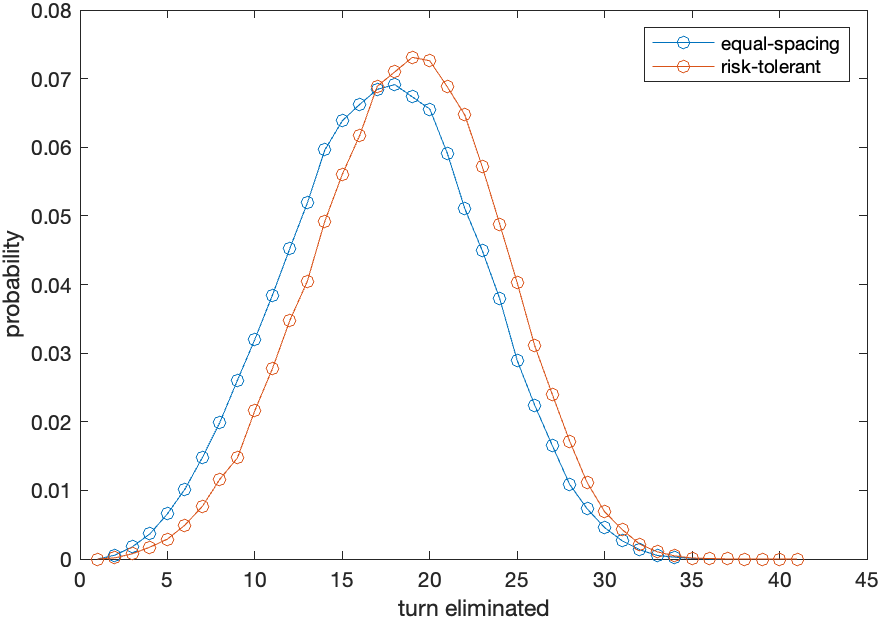}
\caption{Distribution of probability of turn of elimination ($f_n^\sigma$) of an $n$-number game, estimated with 100,000 iterations (\emph{Left}) $n=20$, (\emph{Right}) $n=40$.}
\end{figure}

\section{Conclusion}

The $n$-number game is a probability exercise with an unintuitive, yet provably optimal strategy. While the risk-tolerant strategy improves the odds of winning an individual game, that advantage is tempered when attempting to win a single game as quickly as possible since the improved strategy keeps the player in losing games for longer. The factor of advantage of risk-tolerant over equal-spacing strategy scales approximately linearly.

Several variants of the $n$-number game have emerged that may be intractable to derive analytically optimal strategies for, yet are interesting to inspect nonetheless. One may consider an $n$-number challenge where the random variables are sampled from a non-uniform distribution. A simple strategy to this game may entail converting the distribution to percentiles and following the risk-tolerant strategy on the transformed variables, however since subsets of the distribution are not similar this would require several layers of transforms to execute. Another variant of the $n$-number game that has gained traction recently is the $n^2$-number grid, where instead of ordering random variables in an ordered list, the player places the random variables in an $n\times n$ grid where both rows and columns must be in ascending order. Since there are several ways to arrange the same $n^2$ numbers in an $n\times n$ grid such that rows and columns are sorted, in general the $n^2$-number grid is easier than the $n^2$-number challenge. However, because a grid in progress cannot be recursively decomposed into smaller grid games the way the $n$-number challenge can, it is unclear how to generate an optimal strategy beyond using a greedy metric of maximizing the probability that the remaining numbers satisfy inequalities induced by entries on the grid (see Figure 6).

\begin{figure}
\centering
\begin{TAB}(e,1.25cm,1.25cm){|c:c:c:c:c|}{|c:c:c:c:c|}
 &  &  &  &  \\
130 &  &   &   & 573  \\
  &   &   &   & 761  \\
  &   &   &   &   \\
  &   &   &   &     
\end{TAB}
\hspace{1cm}
\begin{TAB}(e,1.25cm,1.25cm){|c:c:c:c:c|}{|c:c:c:c:c|}
 & 0.16 & 0.34 & 0.41 & 0.35 \\
 & 0.81 & 0.32  & 0.05  &   \\
 0.86 &  0.15 & 0  & 0  &   \\
 0.48 & 0.01  & 0  &  0 &   \\
 0.16 &  0 & 0  &  0 &     
\end{TAB}
\caption{Example of a $5^2$-number grid in progress with 170 as the next number (\emph{Left}) First three turns of the number grid (same entries as in TAB) (\emph{Right}) probability that the remaining 21 numbers satisfy inequalities induced by placing 170 in position $(i,j)$. For example, if 170 is placed in position (3,1), then of the remaining 21 numbers, one must satisfy $x<130$, four must satisfy $x<573$, three must satisfy $130<x<573$, three must satisfy $170<x<761$, eight must satisfy $x>170$, and two must satisfy $x>761$. The probability of this occurring is approximately 0.86 (not to be confused with probability 0.86 of winning from this point).}
\end{figure}

\bibliographystyle{unsrt}
\bibliography{example}

\end{document}